\documentclass[12pt]{article}
\pdfoutput=1  
\usepackage{graphicx}
\usepackage{amssymb, amsmath}
\usepackage{float}

\newtheorem{thm}{Theorem}[section]

\begin{document}
\author{Daniel Carroll, Eleanor Hankins, Emek Kose, Ivan Sterling} 
\title{A Survey of the Differential Geometry of Discrete Curves} 
\date{\today}
\maketitle
\section{Introduction}
Discretization of curves is an ancient topic.  Even discretization of curves with an eye toward differential geometry is over a century old.  However there is no general theory or methodology in the literature, despite the ubiquitous use of discrete curves in mathematics and science.  There are conflicting definitions of even basic concepts such as discrete curvature $\kappa$, discrete torsion $\tau$, or discrete Frenet frame.

Consider for example the three equally worthy definitions of the curvature of an angle derived in Section \ref{CircsNgons} by considering the problem of approximating an N-gon (by N-gon we mean a regular N-gon) with sides of length $\ell$ by a circle: 
\begin{equation} 
\kappa = \frac{2}{\ell} \sin \frac{\theta}{2}, \;\; \kappa = \frac{2}{\ell} \tan \frac{\theta}{2}, \;\;
\kappa = \frac{\theta}{\ell}.
\end{equation}
In the literature each of these definitions occur frequently.  For example \cite{MS}, \cite{H}, \cite{DS}.  As we show, the source of this variety is that each author chooses whether to normalize their curvature by using the inscribed, circumscribed or centered circle of an N-gon.  Although our initial interest was in particular applications, we realized the need for a general approach and along the way discovered some pleasing theorems.

Using Section \ref{CircsNgons} as a guide we proceed to build three theories of discrete curves all of which culminate in a discrete version of the Frenet equations:
\begin{alignat}{4} \label{1disfre}
DT^e &=& \kappa N^v&, \nonumber \\ 
DN^e &=-\kappa T^v &&+\tau B^v, \\
DB^e &=&-\tau N^v&. \nonumber
\end{alignat}
Although dozens of discrete Frenet equations can be found in the literature, all have unpleasant error terms.  Our approach is new, and the resulting Equations (\ref{1disfre}) are free any error terms.  We also show that our definitions of discrete length $\ell$, curvature $\kappa$ and $\tau$ reproduce a unique (up to rigid motion) discrete curve with the given $\ell$, $\kappa$ and $\tau$.


In each of the three cases -- inscribed, circumscribed, and centered -- there corresponds a natural differential geometric way to define the discretization of a smooth curve.  These definitions are discussed in Section \ref{Discretize}.  Conversely given a discrete curve there is a natural differential geometric way to spline the curve.  See Section \ref{Spline}.  Of particular interest is that discrete curves in the plane $\mathbb{R}^2$ are naturally splined by special piecewise curves: constant curvature in the inscribed case, clothoids in the circumscribed case, and elastic curves in the centered case.  In each case we argue for our definition by showing that these splines are the constrained minimizers of the natural variable in that case.  Section \ref{Comments} contains some brief comments about applications and discrete surface theory.

\section{Circles and N-gons} \label{CircsNgons}
\begin{figure}
\begin{center}
\includegraphics[scale=1]{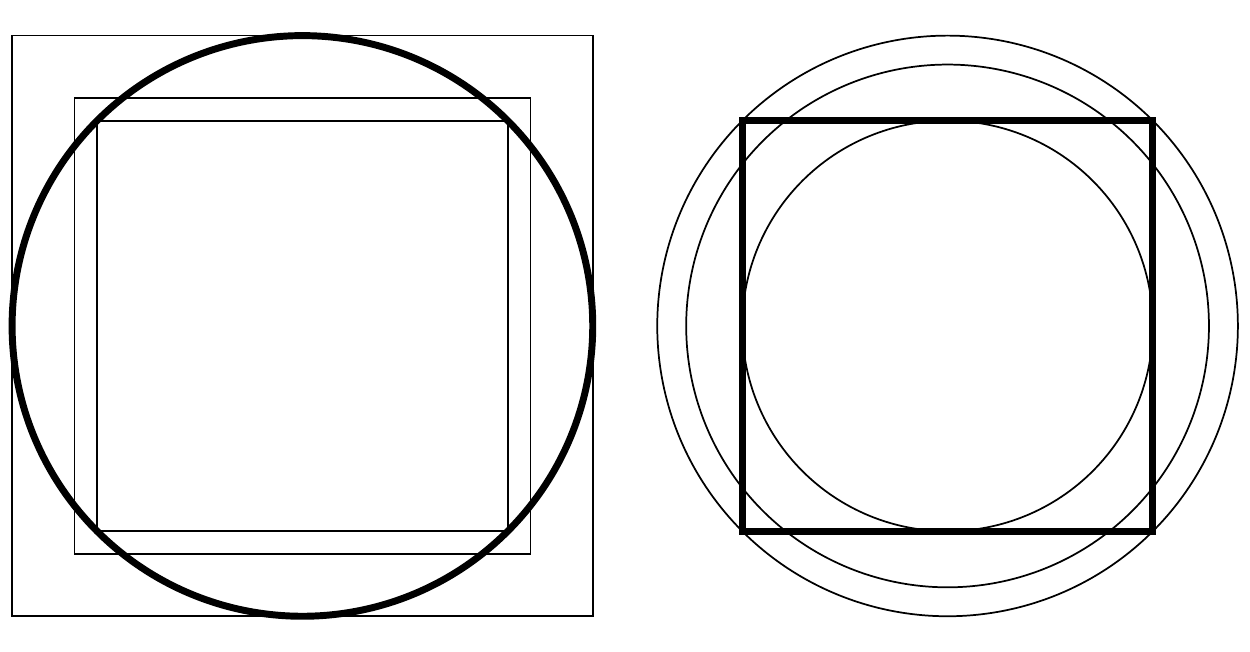}
\caption{Discretizing a Circle, Splining an N-gon}\label{dissplin}
\end{center}
\end{figure}

\begin{figure}
\hspace{-.53in} \includegraphics[scale=.35]{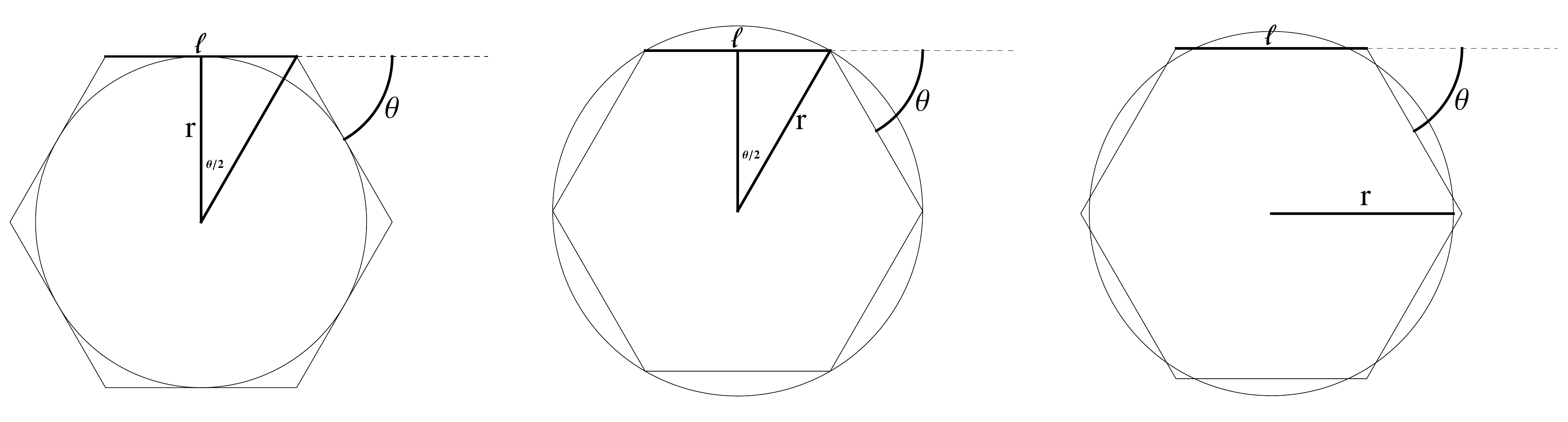}
\caption{Trigonometry for curvature of an N-Gon}\label{trig}
\end{figure}

For every $N \geq 3$, a circle has an inscribed, circumscribed and centered N-gon which discretizes it.  An N-gon is called centered about a circle if its perimeter equals the circumference of the circle.  Conversely an N-gon has an inscribed, a circumscribed and a centered circle which splines it.  See Figure \ref{dissplin} for the case $N=4$.  The nomenclature can be confusing, so care must be taken.  For example: ``An N-gon is centered about a circle" means the circle is given first, whereas ``a circle is inscribed inside an N-gon" means the N-gon is given first.
From trigonometry, Figure \ref{trig}, we are led to three definitions of curvature for a given N-gon.  Recall the curvature of a circle is defined by $\kappa = 1/r$ and that the exterior angle for an N-gon is $\theta = \frac{2 \pi}{N}$.  We assume all sides have length $\ell$ and $N$ is any real $N>0$.  We have the curvature of the circle inscribed in an N-gon is
\begin{equation} \label{k1}
\kappa = \frac{2}{\ell} \sin \frac{\theta}{2}.
\end{equation}
Similarly the curvature of the circle circumscribing an N-gon is
\begin{equation} \label{k2}
\kappa = \frac{2}{\ell} \tan \frac{\theta}{2}, 
\end{equation}
and in the centered case we have
\begin{equation} \label{k3}
\kappa = \frac{\theta}{\ell}.
\end{equation}

\begin{figure}
\hspace{-.2in} \includegraphics[scale=.39]{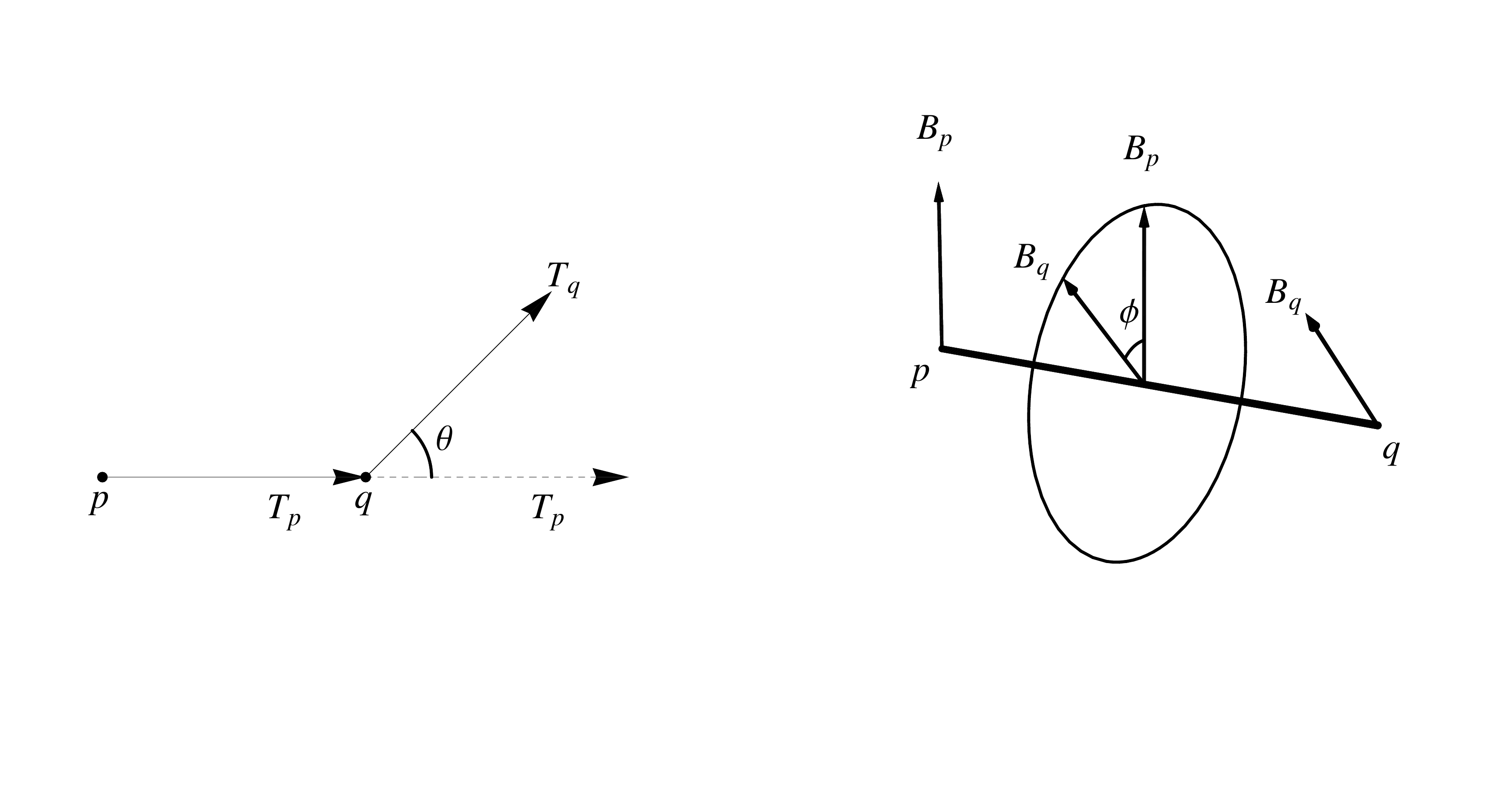}
\vspace{-1in}
\caption{$\theta$ Measures Turning, $\phi$ Measures Twisting}\label{turntwist}
\end{figure}

We will use these basic formulas to guide us in all the definitions that follow. In a way that will be made more precise below we consider $\theta$ as the measure of the angle between neighboring ``tangent vectors".  $\theta$ measures the turning of an N-gon at a vertex.  For a discrete curve in three space if we similarly define $\phi$ to measure the angle between neighboring ``binormal vectors" then $\phi$ measures the twisting of a discrete curve along its edge.  \mbox{See Figures \ref{turntwist} and \ref{dan}.}
We define the curvature at a vertex of a discrete curve in three space by Equations (\ref{k1}), (\ref{k2}), and (\ref{k3}).
We are similarly led to define the torsion at a vertex by
\[\tau := \left\lbrace \begin{array}{l}
 \frac{2}{\ell} \sin \frac{\phi}{2},\;\mbox{in the inscribed case,} \\
 \frac{2}{\ell} \tan \frac{\phi}{2},\;\mbox{in the circumscribing case,} \\
 \frac{\phi}{\ell}, \;\mbox{in the centered case.}
\end{array} \right.\]
 
\section{Discrete Frenet Equations}
A discrete map is a function with domain $\mathbb{Z}$, $\chi: \mathbb{Z} \longrightarrow R$.  Such a map is called a discrete function (resp. curve) if the range is $\mathbb{R}$ (resp. $\mathbb{R}^3$).  Since we work exclusively with these special ranges, we will use without further comment the standard operations of $\mathbb{R}$ and $\mathbb{R}^3$.  If $\chi: \mathbb{Z} \longrightarrow R$, then we often use the notation $\chi_i :=\chi(i)$.  We define discrete differentiation (resp. addition) by $(D\chi)_i := \chi_{i+1}-\chi_i$ (resp. $(M\chi)_i := \chi_{i+1}+\chi_i$).
\subsection{Frenet Frames}
We will define the lengths $\ell_i$, curvatures $\kappa_i$ and torsions $\tau_i$ of discrete curves in such a way that given any $\ell_i$, $\kappa_i$, $\tau_i$ it is possible to reconstruct a discrete curve with these lengths, curvatures and torsions.  We will also require that a natural discrete version of the Frenet equations hold.  As we have seen, there are at least three reasonable definitions of the curvature of the elementary N-gon.  We will investigate these three cases using the definitions of curvature and torsion derived from the formulas above.

\begin{figure}
\hspace{-.2in} \includegraphics[scale=.4]{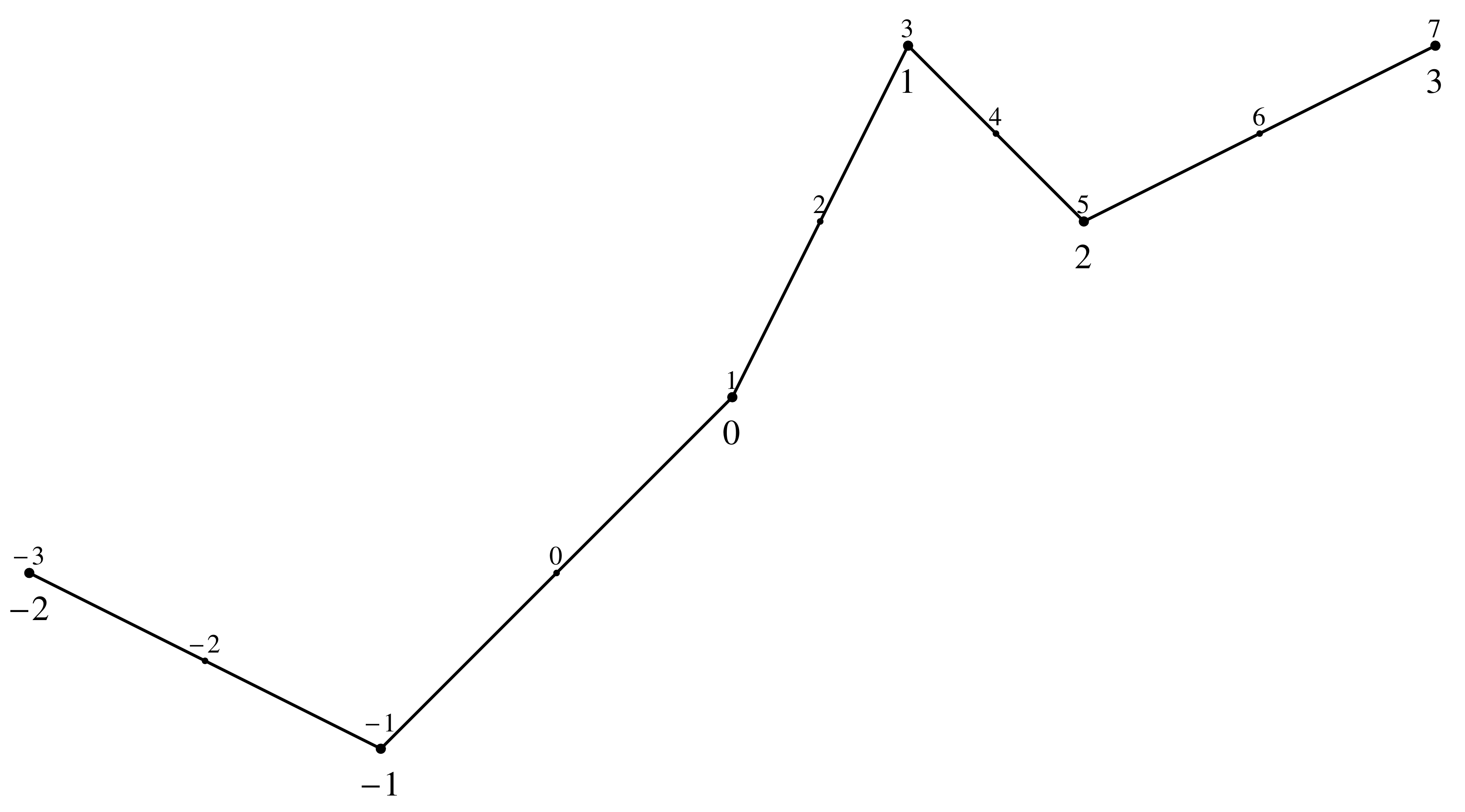}
\vspace{-.25in}
\caption{Original and Redefined Discrete Curve}\label{origredine}
\end{figure}

Let $\gamma^{orig}$ be a discrete curve 
\[\gamma^{orig}: \mathbb{Z} \longrightarrow \mathbb{R}^3,\]
which we call ``the original curve."  Then we define the curve $\gamma:\mathbb{Z} \longrightarrow \mathbb{R}^3$ as follows.  See Figure \ref{origredine} where the larger numbers are the indices for the original curve and the smaller numbers are the indices for the redefined curve.
 First we define
\[\gamma(i):=\gamma^{orig} \left(\frac{i-1}{2}\right) \;\;\; \mbox{\rm{if $i$ is odd}}\]
and then
\[\gamma(i):=\frac{\gamma(i+1) +\gamma(i-1)}{2}  \;\;\; \mbox{\rm{if $i$ is even}}.\]
Note that we recover the original curve from the odd indices of $\gamma$ and that the even indices are mapped to the midpoints of the original curve.  We define the discrete length by 
\[\ell_i := \Vert (D \gamma)_i \Vert.\]
$\gamma$ is parametrized by arc length if $\ell \equiv 1$ and it is parametrized proportional to arc length if $\ell$ is constant.  Note that $\ell_\gamma \equiv \ell = constant$ if $\ell_{\gamma^{orig}} \equiv 2\ell$.  For clarity of presentation we will assume from now on that $\gamma$ is parametrized proportional to arc length,  $\Vert D \gamma \Vert \equiv \ell = constant$.  The theory goes through without this restriction.

\subsection{Frenet Equations}
In each version (Inscribed, Circumscribed and Centered) we will produce two discrete Frenet frames $\{T^e,N^e,B^e\}$ and $\{T^v, N^v,B^v\}$.
First  for $\{T^e,N^e,B^e\}$:
\[T^e:=\frac{D \gamma}{\Vert D \gamma \Vert} = \frac{D \gamma}{\ell}.\]
Note $T^e_i=T^e_{i-1}$ if $i$ is even.
Then 
\[B^e_i := \frac{T^e_i \times T^e_{i+1}}{\Vert T^e_i \times T^e_{i+1} \Vert} \;\;\; \mbox{\rm{if $i$ is even}}\]
and
\[B^e_i := B^e_{i-1} \;\;\; \mbox{\rm{if $i$ is odd}}.\]
and finally for all $i$:
\[N^e_i := B^e_i \times T^e_i.\]
For $\{T^v, N^v,B^v\}$ we have for all $i$:
\[T^v_i:=\frac{(MT^e)_i}{\Vert (MT^e)_i \Vert},\]
\[B^v_i := \frac{(MB^e)_i}{\Vert (MB^e)_i \Vert},\]
\[N^v_i := \frac{(MN^e)_i}{\Vert (MN^e)_i\Vert}.\]
Note for all $i$, $N^v_i = B^v_i \times T^v_i$.

\begin{figure}[H]
\hspace{-1in} \includegraphics[scale=.5]{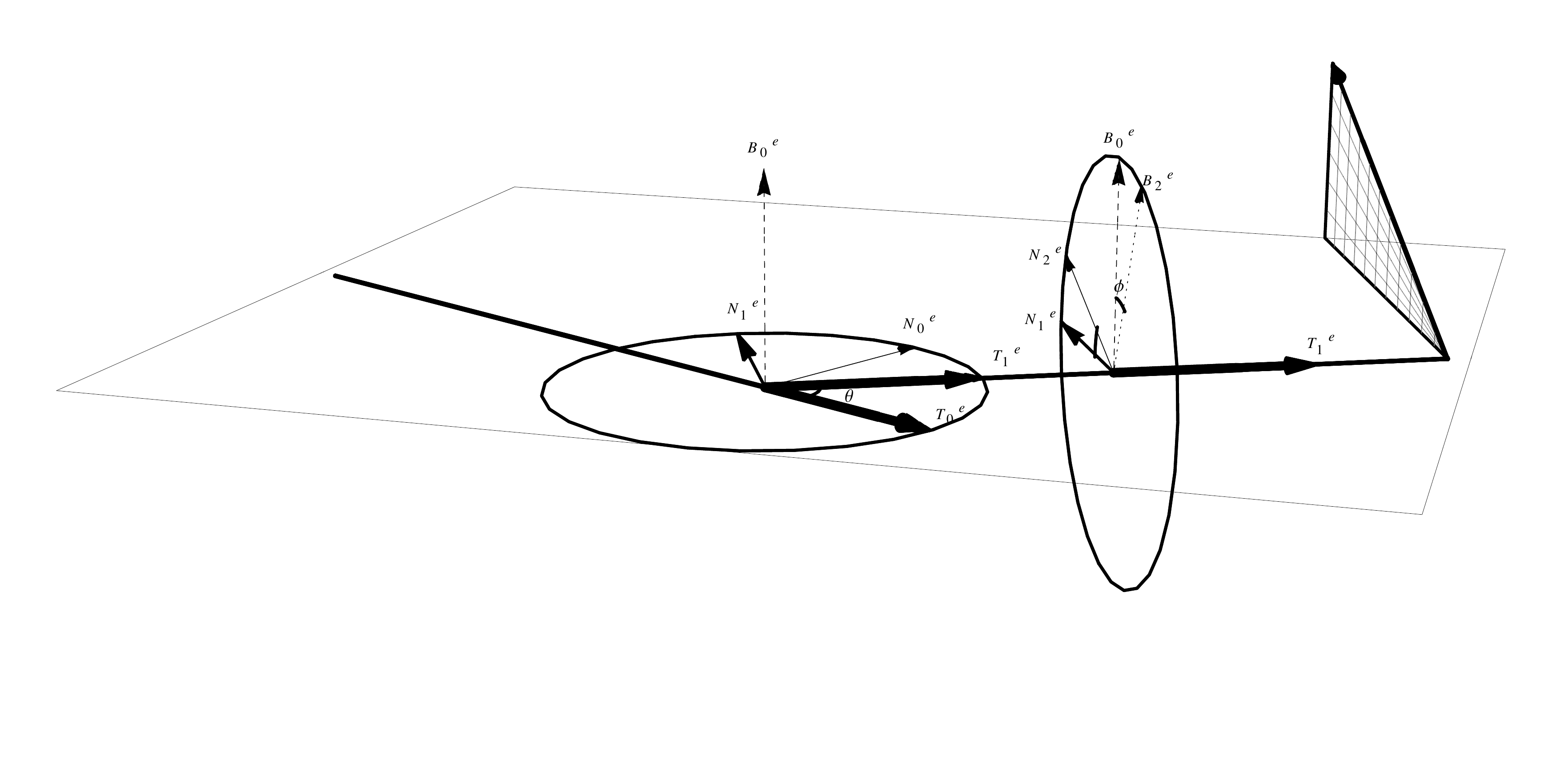}
\vspace{-1in}
\caption{Bird's Eye View}\label{dan}
\end{figure}

As shown again in Figure \ref{dan} the frame is turning, about the axis determined by the binormal, at the ``vertices".  The frame is twisting, about the axis determined by the tangent, at the ``edges".  It was precisely this alternating approach which lead to the elegant form of the discrete Frenet equations (\ref{disfre}) given below; which do not appear in the literature.
\subsection{Curvature and Torsion}
The positively oriented frames {\small $\{T^e_i,N^e_i, B^e_i\}$} determine orientations of {\small $\{T^e_i, N^e_i\}$} and  $\{N^e_i, B^e_i\}$.  We define $\theta_i$ as the angle between $T^e_{i}$ and $T^e_{i+1}$ and note that $\theta_i=0$ if $i$ is odd.  We define $\phi_i$ as the angle between $B^e_{i}$ and $B^e_{i+1}$ with $\phi_i=0$ if $i$ is even.  To avoid technical details we will assume $\theta_i, \phi_i \in [0,\frac{\pi}{2}]$.

The curvature $\kappa$ is defined by
\[\kappa := \left\lbrace \begin{array}{l}
\Vert DT^e \Vert = \frac{2}{\ell} \sin \frac{\theta}{2},\;\mbox{in the inscribed case,} \\
\\
\frac{\|DT^e\|}{\|MT^e\|} = \frac{2}{\ell} \tan \frac{\theta}{2},\;\mbox{in the circumscribing case,} \\
\\
2\sin^{-1}{(\frac{\|DT^e\|}{2})}=\frac{\theta}{\ell}, \;\mbox{in the centered case.}
\end{array} \right.\]
Note that $\kappa_i =0$ if $i$ is odd.  Similarly the torsion $\tau$ is defined by
\[\tau:= \left\lbrace \begin{array}{l}
\phi\; \Vert DB^e \Vert = \frac{2}{\ell} \sin{\frac{\phi}{2}},\;\mbox{in the inscribed case,} \\
\\
 \frac{\|DB^e\|}{\|MB^e\|} = \frac{2}{\ell} \tan \frac{\phi}{2},\;\mbox{in the circumscribing case,} \\
\\
2\sin^{-1}\frac{\|DB^e\|}{2}=\frac{\phi}{\ell}, \;\mbox{in the centered case.}
\end{array} \right.\]
With  $\tau_i=0$ if $i$ is even.
\subsection{Discrete Frenet Equations} \label{secdisfre}
In each version (Inscribed, Circumscribed and Centered) a direct calculation shows that the discrete Frenet equations hold.
\begin{thm}
\begin{alignat}{4} \label{disfre}
DT^e &=& \kappa N^v&,  \nonumber \\
DN^e &=-\kappa T^v &&+\tau B^v, \\
DB^e &=&-\tau N^v&.  \nonumber
\end{alignat}
\end{thm}
\subsection{Discrete Fundamental Theorem}
On the other hand we can reconstruct the curve by the relations:
\begin{alignat*}{4}
T^e_{i+1} &=\;\;\;\cos \theta_i T^e_i+& \sin \theta_i N^e_i & , \\
N^e_{i+1} &=-\sin \theta_i T^e_i+& \sin (\theta_i+\phi_i) N^e_i&- \sin \phi_i B^e_i, \\
B^e_{i+1} &=&\cos \phi_i N^e_i& +\sin \phi_i B^e_i.
\end{alignat*}
and $\gamma_{i+1} = \gamma_i + T^e_{i+1}$.

To summarize we have
\begin{thm} Given $\theta_i$, $\phi_i$ with $\theta_i=0$ for $i$ odd and $\phi_i=0$ for $i$ even.  Then for arbitrary initial conditions $\gamma_0, T^e_0, N^e_0, B^e_0$ there exists a unique discrete curve $\gamma$ with $\theta^\gamma_i = \theta_i, \phi^\gamma_i=\phi_i$ satisfying $\gamma(0)=\gamma_0, {T^\gamma}^e_0= T^e_0, {N^\gamma}^e_0= N^e_0, {B^\gamma}^e_0=B^e_0$.  Moreover, $\gamma^{orig} (i) := \gamma(2 i)$ satisfies $\Vert D \gamma^{orig} \Vert=2 \ell$.
\end{thm}

\section{2D-Discretizing} \label{Discretize}
Given a curve in $\mathbb{R}^2$ we would now like to discretize it.  There is a canonical geometric discretization in each of our three cases.
\subsection{Inscribed 2D-Discretization}  
\begin{figure}[H]\begin{center}\includegraphics[scale=.6]{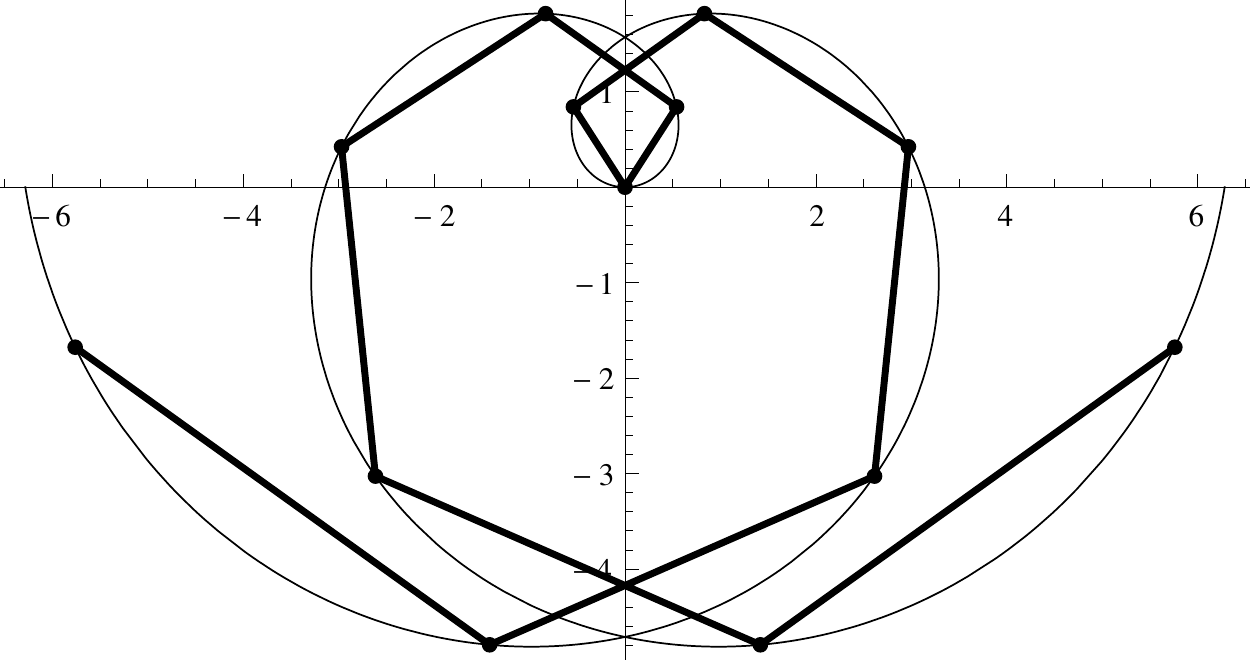}\caption{Inscribed Discretization}\label{ivaninscribe}\end{center}\end{figure}

The only distinguishing feature in this case is that each vertex of the discretization be on the curve itself.  Thus any increasing map $\iota: \mathbb{Z} \longrightarrow \mathbb{R}$ will produce an acceptable discrete curve $\delta := \gamma \circ \iota$.  See Figure \ref{ivaninscribe}.
\subsection{Circumscribed 2D-Discretization}
\begin{figure}[H]\begin{center}\includegraphics[scale=.6]{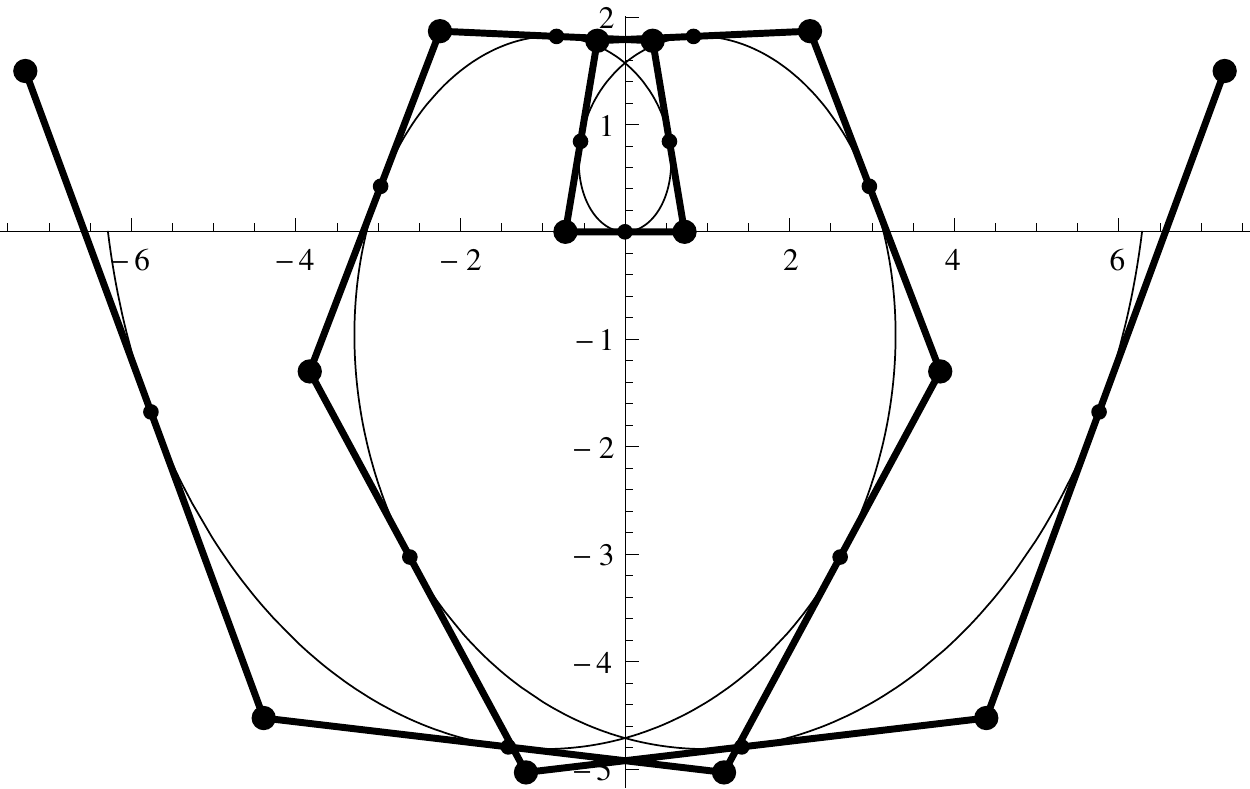}\caption{Circumscribed Discretization}\label{ivancircum}\end{center}\end{figure}

If there are no inflection points then we again take any increasing map $\iota: \mathbb{Z} \longrightarrow \mathbb{R}$ such that consecutive tangents are not parallel.  We require that the edges of our discrete curve $\delta$ intersect tangentially with the given curve at the points $(\gamma \circ \iota)_i$.  We define $\delta_i$ to be the unique intersect point of tangent lines at $(\gamma \circ \iota)_i$ and  $(\gamma \circ \iota)_{i+1}$ as in Figure \ref{ivancircum}.  
\noindent If there are isolated inflection points then they, as well as at least one point in-between them, need to be included in the set of tangent points.
If a curve has infinitely many inflection points on a finite interval, then our algorithm fails.
\subsection{Centered 2D-Discretization}
\begin{figure}[H]\begin{center}\includegraphics[scale=.6]{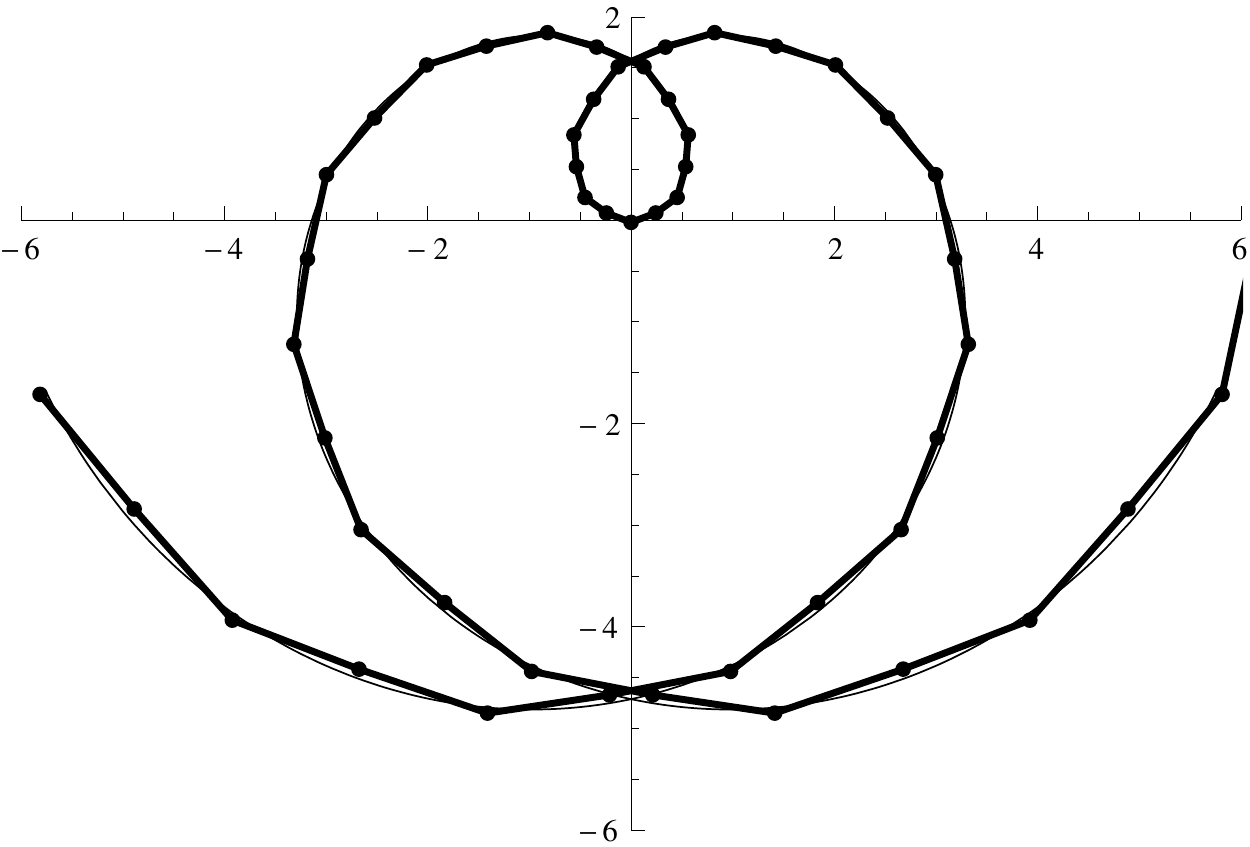}\caption{Centered Discretization}\label{perfected}\end{center}\end{figure}

\begin{figure}[H]\begin{center}\includegraphics[scale=.6]{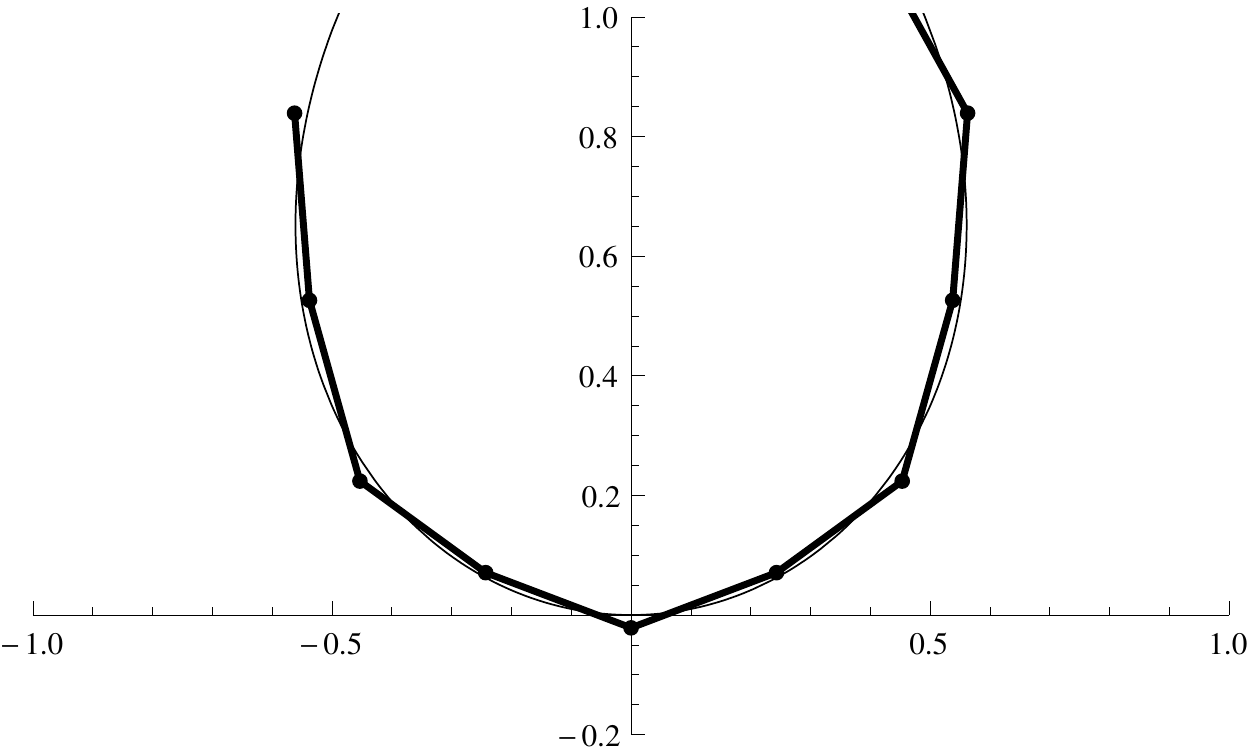}\caption{Centered Discretization Zoom}\label{perfectedtiny}\end{center}\end{figure}

The natural centered discretization of a curve requires a bit more finesse.  First, without loss of generality, we assume $\gamma$ is parametrized by arc length and require our discretization to be parametrized proportional to arc length.  Secondly, with loss of generality, we assume $\gamma$ has no inflection points, say $\kappa > 0$ everywhere.  We require $\kappa_i > 0$ for our discretization.  Finally we choose $M$ ``large enough."  We take the specific $\iota: \mathbb{Z} \longrightarrow \mathbb{R}$ defined by $\iota(i):=\frac{i}{M}$ and let $\delta^{start} = \gamma \circ \iota$.  For each $i$ we offset $\delta^{start}_i$ along the (outward) normal to $\gamma$ at $\delta^{start}_i$ by the amount \[\mbox{offset}_i := \frac{\frac{k_i}{M} - \sin{\frac{k_i}{M}}}{k_i \sin{\frac{k_i}{M}}},\]
where $k_i$ is the curvature of $\gamma$ at $\delta^{start}_i$ and we assume $k_i > 0$.  Note that this formula is derived from the case of centered N-gons discretizing a circle.
\begin{figure}[H]\begin{center}\includegraphics[scale=.6]{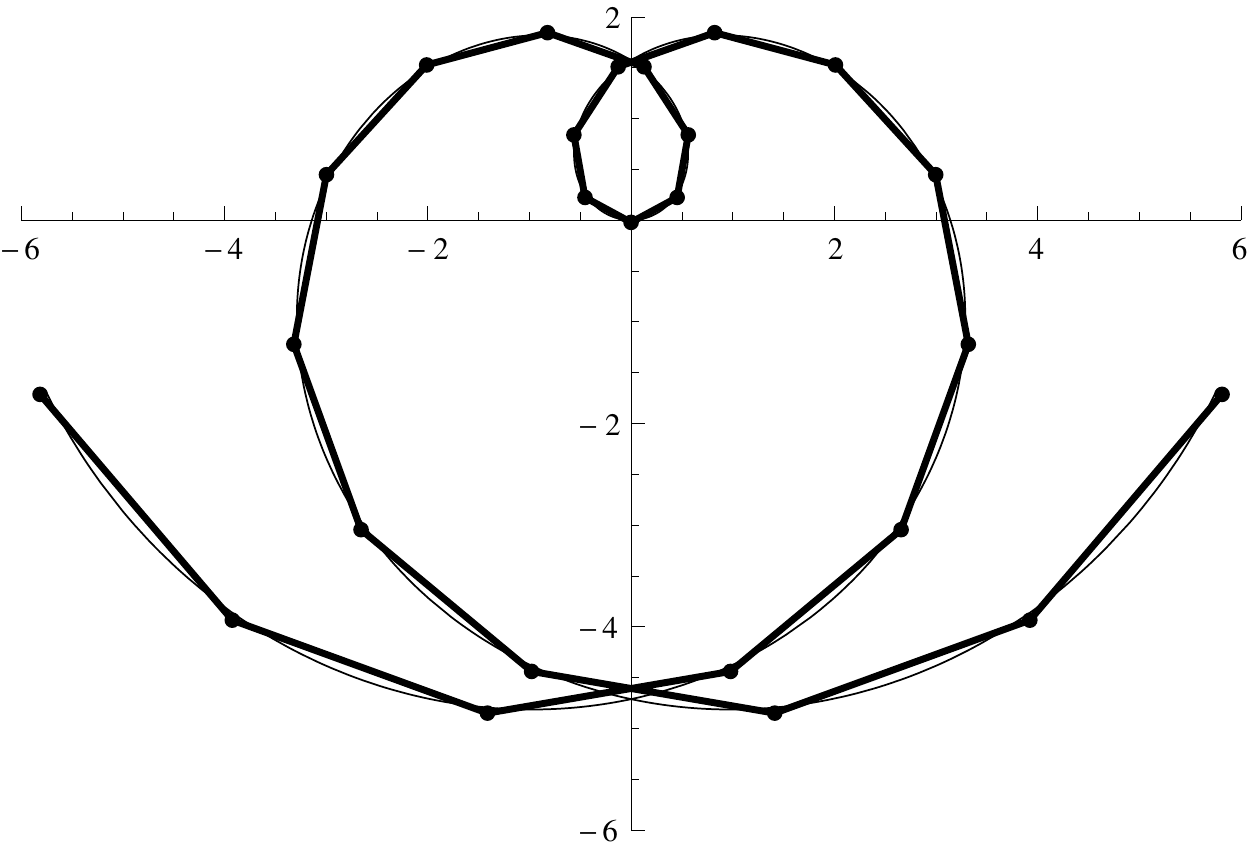}\caption{Offset Discretization}\label{offsetonly}\end{center}\end{figure} 
\noindent To see the offset more clearly we zoom into the center of the curve.
\begin{figure}[H]\begin{center}\includegraphics[scale=.6]{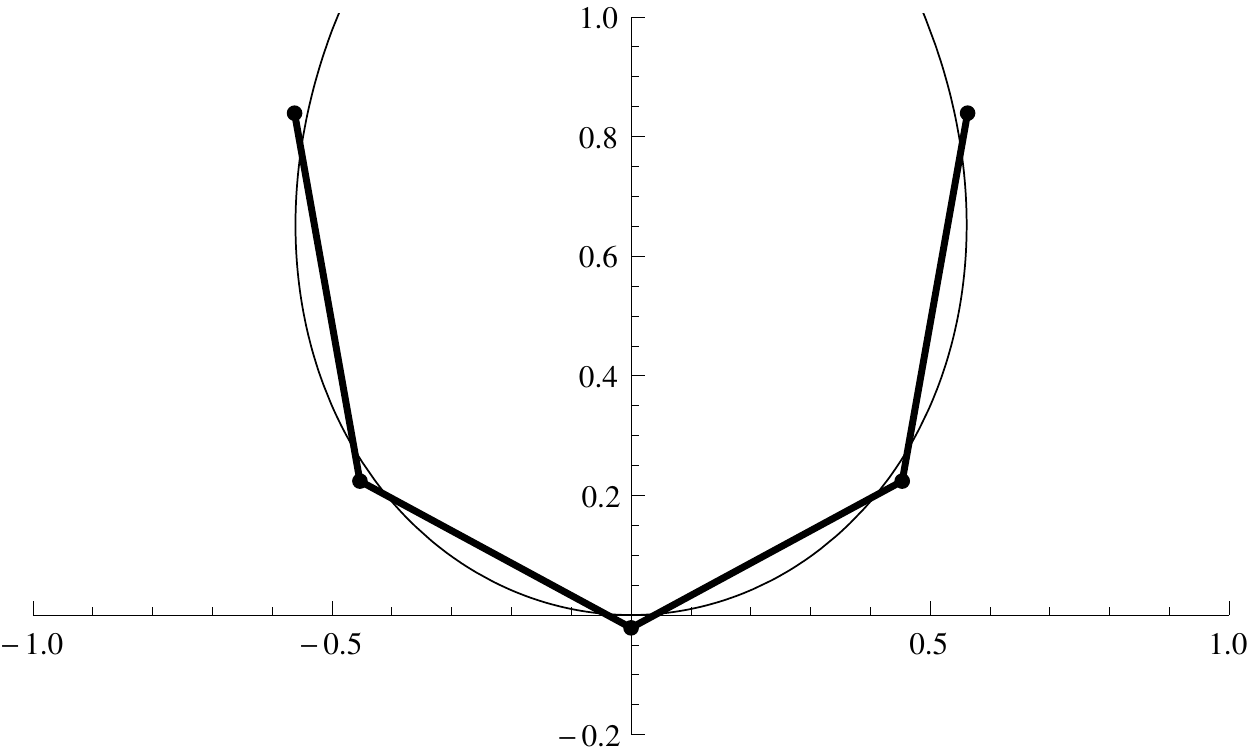}\caption{Offset Discretization Zoom}\label{offsetonlytiny}\end{center}\end{figure}
\noindent These offset points will be the even vertices, $\delta_{2j}$, of our final discrete curve $\delta$.   Now we consider the condition that our discrete curve $\delta$ is to have the same length as our original curve $\gamma$.  With the additional conditions that  $\Vert \delta_{2j+1}-\delta_{2j} \Vert=\frac{i}{2M}$, $\Vert \delta_{2j+2}-\delta_{2j+1} \Vert$ and $\kappa_{2j+1} > 0$; we see there is one and only one way to achieve this.   See Figure  \ref{perfected}.  
\noindent Again, we see more detail by zooming in, Figure \ref{perfectedtiny}.
\noindent To include inflection points requires more general parametrizations and we will leave it as an exercise for the reader.

\section{Geometric Splinings of Discrete Curves} \label{Spline}
\begin{figure}[H]\begin{center}\includegraphics[scale=0.5]{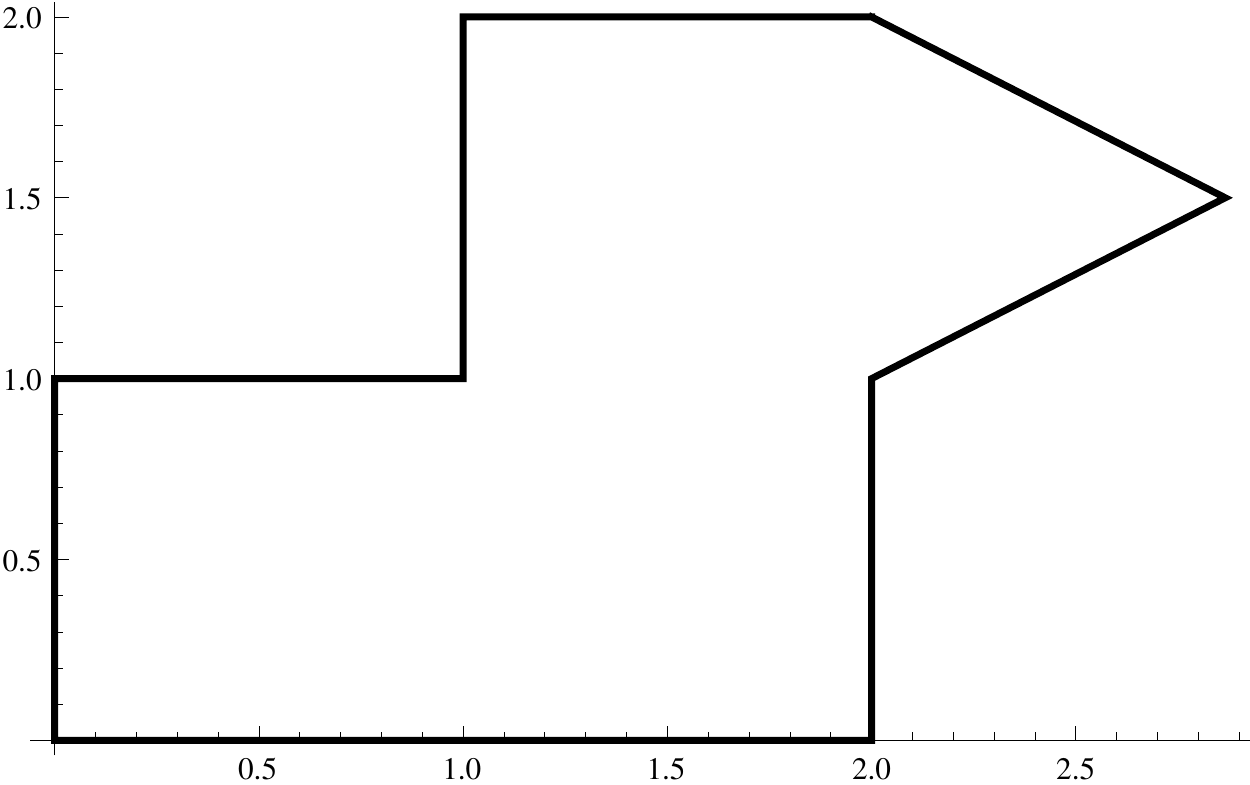}
\caption{Discrete Curve to be Splined}\label{emekpreinscribed}\end{center}\end{figure}
\subsection{What does Best Spline Mean?}
In non-geometric splining, the best spline is usually related to the degree of the polynomial $\gamma(t)=(x(t),y(t))$ used to approximate the curve.  For example a cubic spline is constructed using piece-wise cubic polynomials.  Typically a cubic spline passes through the points of a discrete curve with certain boundary conditions.  Geometric splinings on the other hand are found by considering curves whose curvature function $\kappa(t)$ is a low degree polynomial.  Alternatively a best geometric spline minimizes $\int \kappa^2$.
\subsection{Inscribed Splining} 
\begin{figure}[H]
\begin{center}
\includegraphics[scale=0.5]{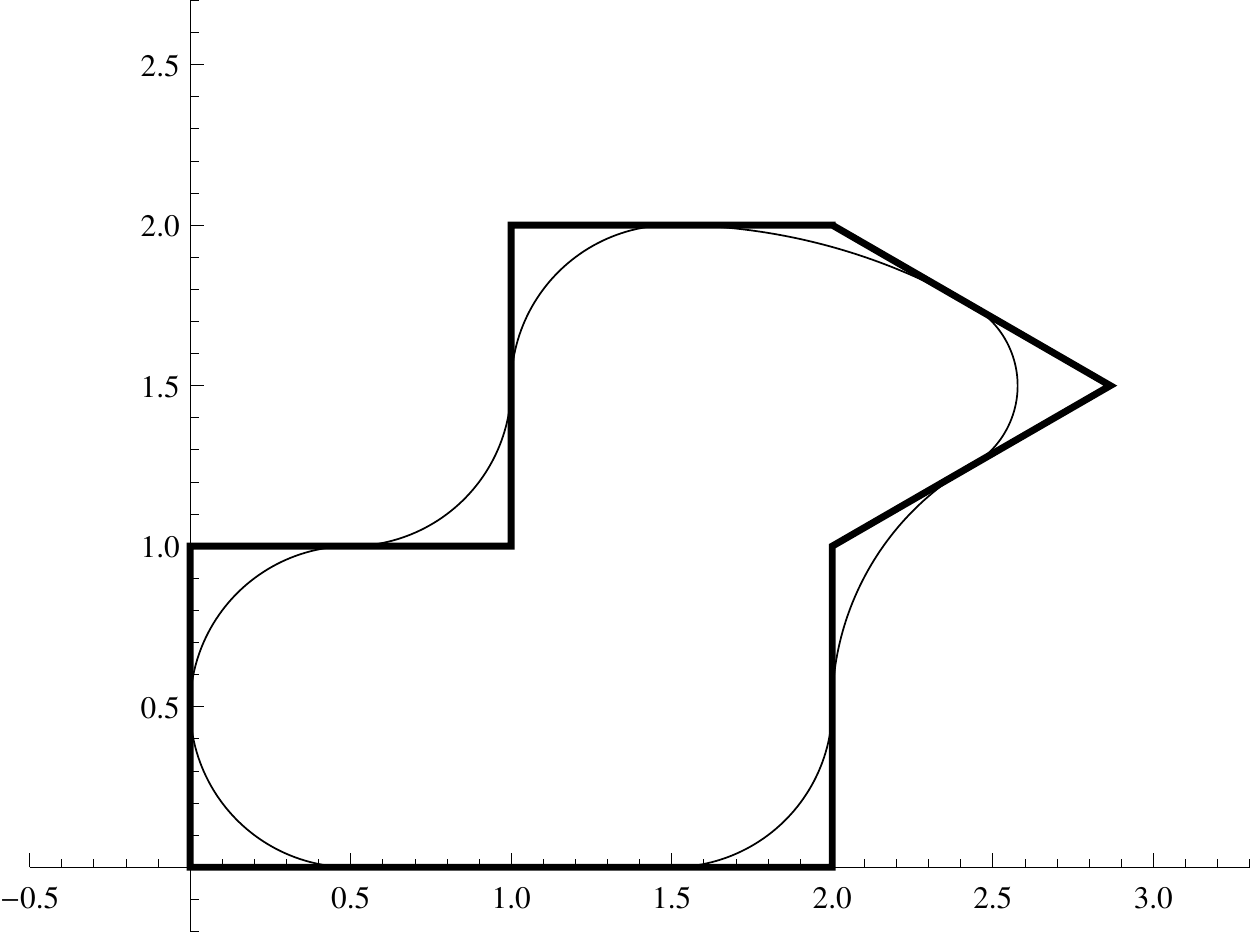}
\caption{Inscribed Splining}\label{newinscribed}
\end{center}
\end{figure}

An inscribing spline is one which tangentially goes through the midpoints of the edges of the given discrete curve.  We seek a curve whose curvature has the lowest degree possible.  Because we are assuming our discrete curves are parametized proportional to arc length there is a trivial differentiable inscribed splining by pieces of curves of constant curvature.  That is pieces of circles.  See Figure \ref{newinscribed}.  If our discrete curve is not parametrized proportional to arc length, then the inscribed splining would require clothoids, which are described in the next subsection.
\noindent Notice the curvature jumps at the midpoints so our splining is not twice differentiable.
\subsection{Clothoids}
Curves with linear curvature are called first order clothoids.  Given
$$\kappa (s) = as+b$$
(if $a=0$,  we get a piece of a circle, a ``zeroth order clothoid") then the turning angle $\theta$ is given by
\[\theta (s) = \int_0^s \kappa(t) \, dt +  \theta_0.\]
First order clothoids are given in terms of Fresnel integrals
\[\gamma(s) = \left(\displaystyle\int_0^s \, \cos{\theta(t)} \,dt + x_0, \displaystyle\int_0^s \, \sin{\theta(t)} \,dt +y_0 \right).\]
Similarly curves with quadratic curvature are second order clothoids, and so on.

\subsection{Circumscribed Splining}
\begin{figure}[H]
\begin{center}
\includegraphics[scale=0.5]{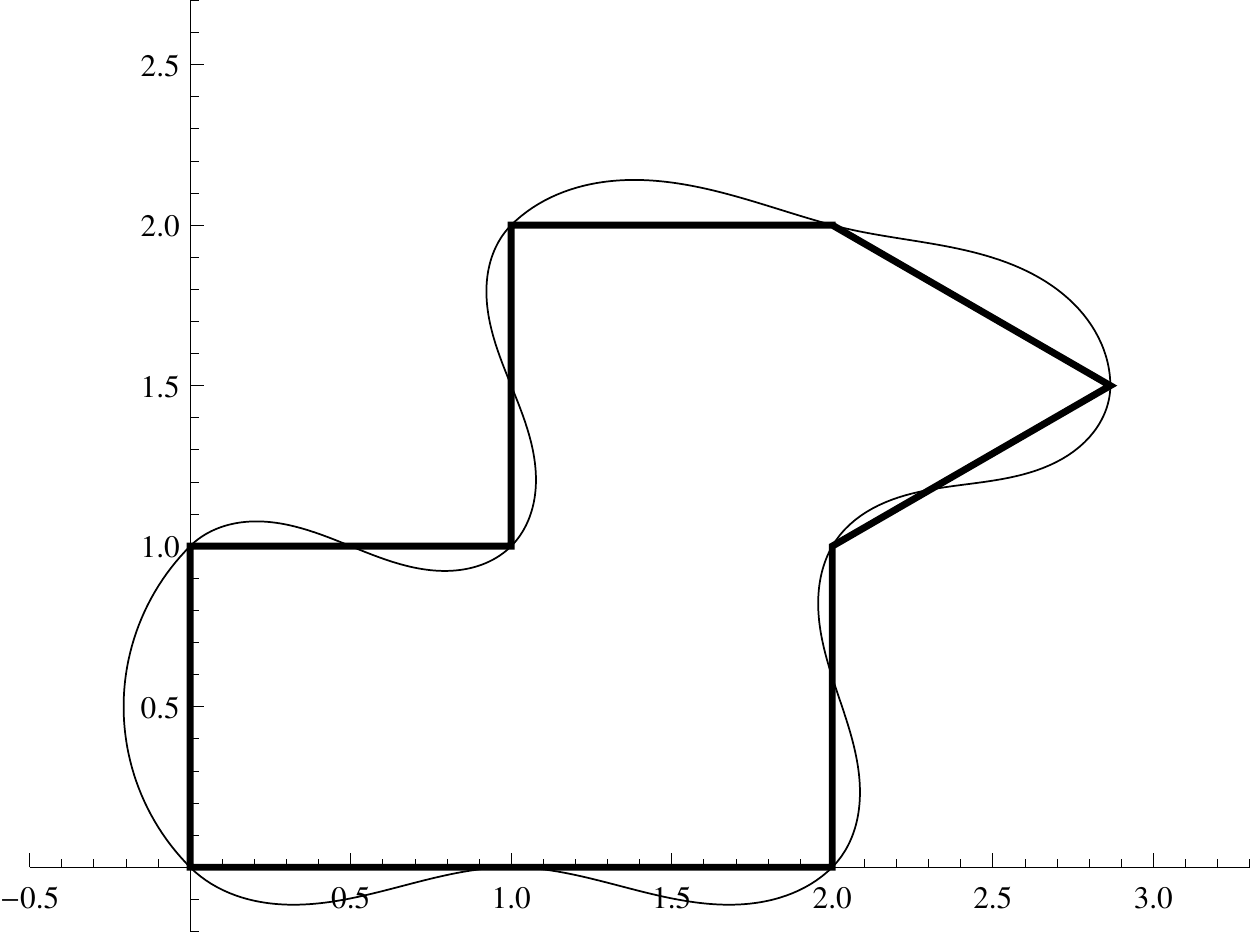}
\caption{Circumscribed Splining}\label{newcircum}
\end{center}
\end{figure}

A circumscribing spline is one which differentiable goes through the points of the given discrete curve.  Unlike the case of inscribed splinings, it will rarely be the case that a circumscribed splining will consist of piece of circles.  On the other hand there will always be circumscribing splining, as in Figure \ref{newcircum}, using first order clothoids.  If there is more than one, we take the shortest one.  This is called the fitting problem.  See for example \cite{BF}.
\noindent Notice again the curvature jumps at the midpoints.
\subsection{Centered Splining}
\begin{figure}[H]
\begin{center}
\includegraphics[scale=0.5]{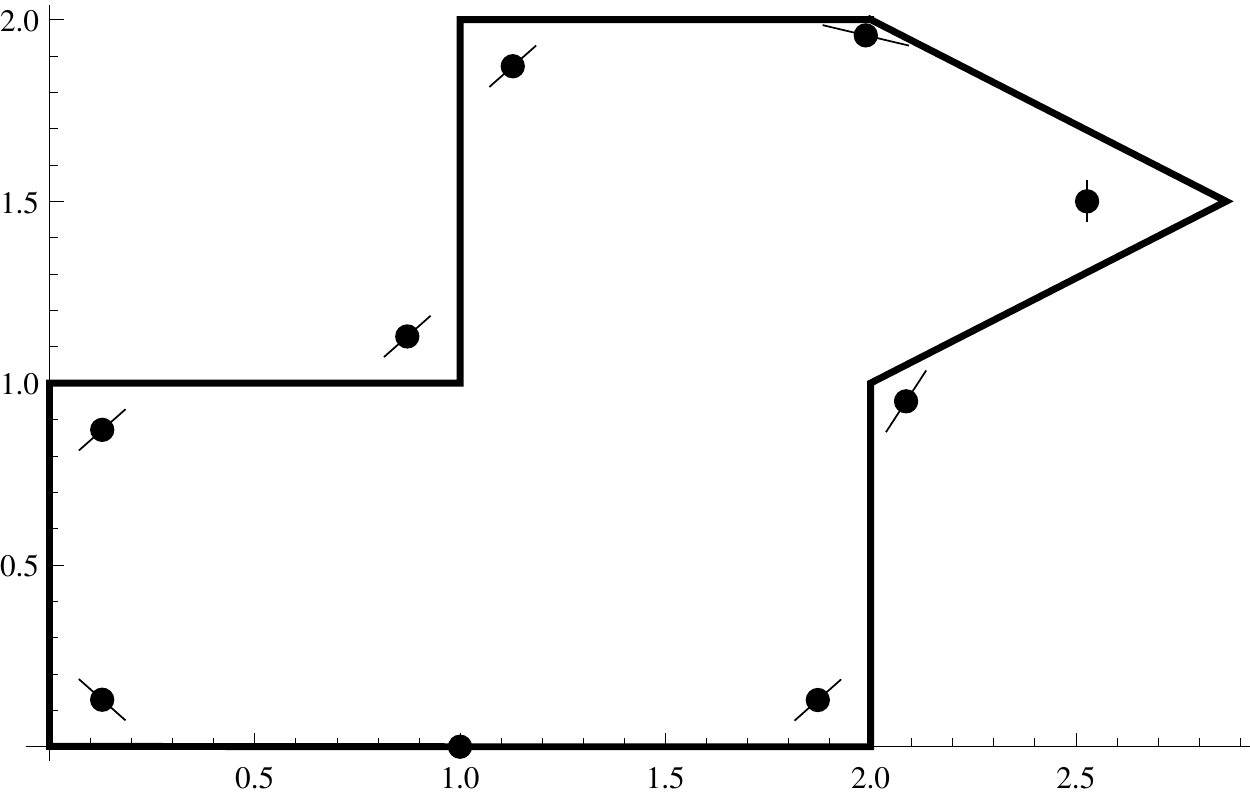}
\caption{Centered Splining Offsets}\label{newcentoffset}
\end{center}
\end{figure}

For the centered spline we first offset the vertices using the centered circles of N-gons and take the directions of the desired spline at these offset point to be the average of the incoming and outgoing directions of the edges at the vertices.  See Figure \ref{newcentoffset}.

We then seek differentiable splines passing through these offset points whose length agrees with that of the given discrete curve.  These curves are found using $\int \kappa^2$ and are called elastica.  These are solutions to a variational problem proposed by Bernoulli to Euler in $1744$; that of minimizing the bending energy of a thin inextensible wire.  Among all curves of the same length that not only pass through points $A$ and $B$ but are also tangent to given straight lines at these points, it is defined as the one minimizing the value of the expression $\int \kappa^2$. 

The one parameter family of elastic curves introduced by Euler \cite{O} is well known.  They are all given by explicit formulas involving elliptic integrals.  These formulas arise by solving the one-dimensional sine-Gordon differential equation.  Which is alternatively written as $\theta''=\sin{\theta}$ or \cite{L} $\theta'''+\frac{1}{2}(\theta')^3 + C \theta' =0$.  In applied problems, such as finding the elastic curve with the boundary conditions care must be taken for several reasons.  One issue is that there are several types of  ``elastic intervals" (inflectional, non-inflectional, critical, circular, and linear).  Another issue is that in some cases there are multiple solutions.  An excellent survey of the subject is in Andentov \cite{A1}.  As discussed in \cite{A1} these problems persist when attempting to numerically approximate elastic curves.

Sogo \cite{S} shows how, at least in some cases, ``integrable discretization" theory can be used to construct a discretized one-dimensional sine-Gordon equation satisfied by discretized elliptic integrals.  For example the inflectional type elastic curve has a turning angle which is given by a formula, involving the Jacobi sn function, of the form
\[\sin \frac{\theta}{2} = \sin \frac{\theta_0}{2} \mbox{sn}(\frac{K}{L}(L-s),k) \]
and Sogo shows that 
\[\sin \frac{\theta_j}{2} = \sin \frac{\theta_0}{2} \mbox{sn}(\frac{K}{N}(N-j),k) \]
is the turning angle of an approximating discrete elastic curve.  

\begin{figure}[H]
\begin{center}
\includegraphics[scale=0.5]{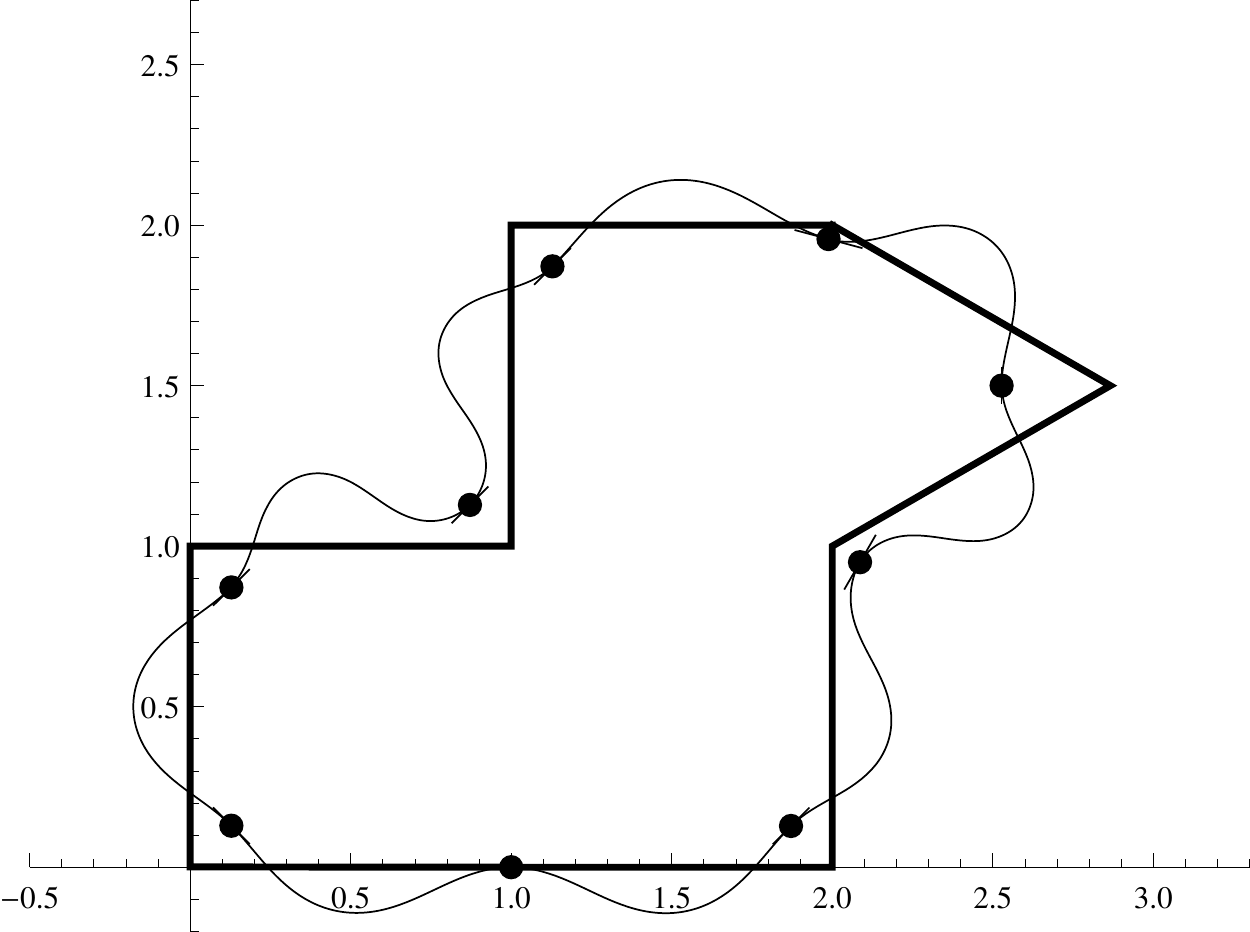}
\caption{Centered Splining}\label{newelastic}
\end{center}
\end{figure}

Figure \ref{newelastic} shows (one of) the differentiable elastic splines with minimal bending energy and length nine.

\section{Comments} \label{Comments}
\begin{itemize}
\item The discretization of smooth curves and the splining of discrete curves in three space is also well studied.  Similar, though at times more involved, case by case constructions can be carried out in all six cases discussed in detail for the curves in the plane.  A complete understanding of circles and N-gons guides the discretization and splining methods in the plane.  Similarly by first carrying out the most basic case of the helix one is able to succeed with curves in three space  as well.
\item The three settings (Inscribed, Circumscribed, and Centered) and only these three settings, are used extensively in the literature.  This is true in both pure and applied differential geometry.  Other settings we considered, although formally feasible, are not as natural.  For example the reader can consider circles whose enclosed areas agree with the enclosed areas of a regular $N$-gons.
\item We feel there is no absolute ``right" definition of discrete curvature or torsion.  A particular application may inform the researcher as to which definitions to use.  For example clothoids arise in the building of highway off ramps.  So in that case the circumscribed setting might be more natural. In a more abstract context the circumscribed setting is also used by T. Hoffman in his dissertation on discrete curves and surfaces \cite{H}.  It seems clear that Gauss would have used the centered setting, as it agrees most closely with his definition of the curvature of an angle between two intersecting curves and with his definition of curvature given by the normal Gauss map.  This setting is used, for example, by Doliwa and Santini \cite{DS} in their work on the integrable dynamics of discrete curves.
\item There is also a vast literature on discrete surface theory which goes back over one hundred years.  See \cite{BP} and references there.  Not surprisingly there is an even wider variety of definitions for the standard concepts such as discrete Gauss curvature, discrete mean curvature, discrete umbilics, etc.  Again it seems clear that there is no absolute ``right" definition.  How one chooses to define ``the discretization" of a smooth surface will again depend on which properties one wishes to preserve.  The theory of ``integrable discretizations" in particular has been applied to soap bubbles, minimal surfaces, Hasimoto surfaces (i.e. the surfaces swept out by smoke-rings) and surfaces of constant Gauss curvature.  Similar comments apply to the theory of splining discrete surfaces.
\item We have highlighted the Frenet frame because it is the most well known curve framing.  Discrete versions of the Bishop frame \cite{B},\cite{CKS} can also be derived using the ideas of this paper.  The Bishop frame is particularly useful for curves that have points of zero curvature.
\item We have considered only the simplest discretizations and the simplest splinings.  One which are as local as possible, taking into account only the ``nearest neighbors."   We feel the diversity and elegance of the cases covered give a nice survey.  Third order versions, either taking into account more points for each calculation or by including curvature into boundary conditions, can be found in the literature.
\end{itemize}


\begin{thebibliography}{99}
\bibitem{A1} Ardentov, A. A. and Sachkov, Yu. L. {\em Solution to Euler’s Elastic Problem}, Automation and Remote Control, Vol 70, No 4 (2009)  633-643.
\bibitem{B} Bishop, R. {\em There is more than one way to frame a curve}, Amer. Math. Monthly 82 (1975) 246-251.
\bibitem{BF} Bertolazzi, E. and  M. Frego, M. {\em Fast and accurate $G^1$ fitting of clothoid curves}, arXiv:1305.6644v2 (2013).
\bibitem{BP} A. Bobenko, A. and  Pinkall, U. {\em Discrete surfaces with constant negative Gaussian curvature and
the Hirota equation}, J. Diff. Geom. 43 (1990) 527-611.
\bibitem{CKS} Carroll, D., Kose, E. and Sterling, I. {\em Improving Frenet's Frame Using Bishop's Frame}, J. of Mathematics Research, Vol 5, No 4 (2013) 97-106.
\bibitem{DS} Doliwa, A. and Santini P. {\em Integrable dynamics of a discrete curve and the Ablowitz-Ladik hierarchy}, J. Math. Phys. 36  (1995) 1259-1273.
\bibitem{H} Hoffman, T. {\em Discrete Curves and Surfaces}, PhD Thesis, Technische Universit\"{a}t Berlin (2000).
\bibitem{L} Levien, R. {\em From Spiral to Spline: Optimal Techniques in Interactive Curve Design}, Publisher, University of California, Berkeley (2009).
\bibitem{MS} McCrae, J. and Singh, K. {\em Sketching Piecewise Clothoid Curves}, Computers \& Graphics, Vol 33, Issue 4 (Aug 2009) 452–461.
\bibitem{O} Oldfather, W. A., Ellis, C. A. and Brown, Donald M. {\em Leonhard Euler's Elastic Curves},  Isis, Vol 20, No 1 (Nov 1933) 72-160.
\bibitem{S} Sogo, K. {\em  Variational Discretization of Euler’s Elastica Problem}, J. of the Physical Society of Japan, Vol 75, Issue 6 (2006) 064007-064007-4.
\end{thebibliography}
\end{document}